\documentclass[12pt]{article}
\usepackage{amssymb}
\newcommand{\be}{\begin{equation}}
\newcommand{\ee}{\end{equation}}
\newcommand{\bea}{\begin{eqnarray}}
\newcommand{\eea}{\end{eqnarray}}
\newcommand{\ba}{\begin{array}}
\newcommand{\ea}{\end{array}}

\newcommand{\bc}{\begin{center}}
\newcommand{\ec}{\end{center}}
\newcommand{\ben}{\begin{enumerate}}
\newcommand{\een}{\end{enumerate}}
\newcommand{\bfi}{\begin{figure}}
\newcommand{\efi}{\end{figure}}

\newcommand{\bq}{\begin{quote}}
\newcommand{\eq}{\end{quote}}
\newcommand{\bqu}{\begin{quotation}}
\newcommand{\equ}{\end{quotation}}
\newenvironment{emphit}{\begin{itemize}}{\end{itemize}}
\newcommand{\bemp}{\begin{emphit}}
\newcommand{\eemp}{\end{emphit}}

\newcommand{\bt}{\begin{tabular}}
\newcommand{\et}{\end{tabular}}

\newtheorem{myth}{Theorem}[section]
\newtheorem{mylem}{Lemma}[section]
\newtheorem{mycor}{Corollary}[section]

\newtheorem{mydef}{Definition}[section]

\begin{document}
\date{}
\title{On Isotopic Characterization of Central Loops
}
\author{T. G. Jaiy\'e\d ol\'a\thanks{All correspondence to be addressed to this author.jaiyeolatemitope@yahoo.com,~tjayeola@oauife.edu.ng} \\
Department of Mathematics,\\
Obafemi Awolowo University,\\
Ile Ife, Nigeria. \and
J. O. Ad\'en\'iran\thanks{E-mail:~ekenedilichineke@yahoo.com} \\
Department of Mathematics,\\
University of Abeokuta, \\
Abeokuta 110101, Nigeria.} \maketitle
\begin{abstract}
The representation sets of central loops are investigated and the
results obtained are used to construct a finite C-loop. It is shown
that for certain types of isotopisms, the central identities are
isotopic invariant.
\end{abstract}

\section{Introduction}
\paragraph{}The isotopic invariance or the universality of types and
varieties of quasigroups and loops described by one or more
equivalent identities, especially those that fall in the class of
Bol-Moufang type loops as first named by Fenyves \cite{phd56} and
\cite{phd50} in the 1960s  and later on in this $21^\textrm{st}$
century by Phillips and Vojt\v echovsk\'y \cite{phd9}, \cite{phd61}
and \cite{phd124} have been of interest to researchers in loop
theory in the recent past. Falconer \cite{phd159} and \cite{phd160}
investigated isotopy invariants in quasigroups and loops. Loops such
as Bol loops, Moufang loops, extra loops are the most popular loops
of Bol-Moufang type whose isotopic invariance have been considered.
But for LC-loops, RC-loops and C-loops, up till this moment, there
is no outstanding result on their isotopic invariance.

Bol-Moufang type of quasigroups(loops) are not the only
quasigroups(loops) whose universality have been considered. Some
others are flexible loops, F-quasigroups, totally symmetric
quasigroups, distributive quasigroups, weak inverse property
loops(WIPLs), cross inverse property loops(CIPLs), semi-automorphic
inverse property loops(SAIPLs) and inverse property loops(IPLs). As
shown in Bruck \cite{phd3}, a left(right) inverse property loop is
universal if and only if it is a left(right) Bol loop, so an IPL is
universal if and only if it is a Moufang loop. Recently, Kepka et.
al. \cite{phd95}, \cite{phd118}, \cite{phd119} solved the Belousov
problem concerning the universality of F-quasigroup which has been
open since 1967. The universality of WIPLs and CIPLs have been
addressed by OSborn \cite{phd89} and Artzy \cite{phd30} respectively
while the universality of elasticity(flexibility) was studied by
Syrbu \cite{phd97}. Basarab \cite{phd146} later continued the work
of Osborn on universal WIPLs by studying isotopes of WIPLs that are
also WIPLs after he had studied a class of WIPLs(\cite{phd149}). The
universality of SAIPLs and Osborn loops is still an open problem to
be solved as stated by Kinyon during the LOOPs'99 conference and
Milehigh conference 2005(\cite{phd33}) respectively. After the
consideration of universal AIPLs by Karklinsh and Klin
\cite{phd150}, Basarab \cite{phd147} obtained a sufficient condition
for which a universal AIPL is a G-loop.

\paragraph{}
LC-loops, RC-loops and C-loops are loops that satisfy the
identities
\begin{displaymath}
(xx)(yz)=(x(xy))z,~(zy)(xx)=z((yx)x)~\textrm{and}~x(y(yz))=((xy)y)z~\textrm{respectively}.
\end{displaymath}
These three types of loops shall be collectively called central
loops. In the theory of loops, central loops are some of the least
studied loops. They have been studied by Phillips and Vojt\v
echovsk\'y \cite{phd9}, \cite{phd61}, \cite{phd58}, Kinyon et. al.
\cite{phd59}, \cite{phd22}, \cite{phd124}, Ramamurthi and Solarin
\cite{ram}, Fenyves \cite{phd56} and Beg \cite{phd169},
\cite{phd170}. The difficulty in studying them is as a result of the
nature of the identities defining them when compared with other
Bol-Moufang identities. It can be noticed that in the aforementioned
LC identity, the two $x$ variables are consecutively positioned and
neither $y$ nor $z$ is between them. A similarly observation is true
in the other two identities(i.e the RC and C identities). But this
observation is not true in the identities defining Bol loops,
Moufang loops and extra loops. Fenyves \cite{phd56} gave three
equivalent identities that define LC-loops, three equivalent
identities that define RC-loops and only one identity that defines
C-loops. But recently, Phillips and Vojt\v echovsk\'y \cite{phd9},
\cite{phd61} gave four equivalent identities that define LC-loops
and four equivalent identities that define RC-loops. Three of the
four identities given by Phillips and Vojt\v echovsk\'y are the same
as the three already given by Fenyves.

Their basic properties are found in \cite{phd58}, \cite{ram},
\cite{phd56} and \cite{den}. The left(right) representation of a
loop $L$ denoted by $\Pi_\lambda(\Pi_\rho)$ is the set of all
left(right) translation maps on the loop i.e if $L$ is a loop, then
\begin{displaymath}
\Pi_\lambda=\{L_x~:~L\to L~\vert~ x\in
L\}~\textrm{and}~\Pi_\rho=\{R_x~:~L\to L~\vert~x\in L\}
\end{displaymath}
where $R_x:L\to L$ and $L_x:L\to L$, defined as $yR_x=yx$ and
$yL_x=xy$ respectively $\forall~x,y\in L$ are bijections. $L$ is
said to be a central square loop if~$x^2\in Z(L, \cdot )~\forall
~x\in L$ where $Z(L,\cdot )$ is the center of $L$. $L$ is said to be
left alternative if~ $\forall~ x, y\in L,~ x\cdot xy=x^2y$ and is
said to be right alternative if~ $\forall~ x, y\in L,~ yx\cdot
x=yx^2$. Thus, $L$ is said to be alternative if it is both left and
right alternative. The set $S(L, \cdot )$ of all bijections in a
loop $(L,\cdot )$ forms a group called the permutation group of the
loop $(L,\cdot )$. The triple $(U, V, W)$ such that $U, V, W\in S(L,
\cdot )$ is called an autotopism of $L$ if and only if
\begin{displaymath}
xU\cdot yV=(x\cdot y)W~\forall~x,y\in L.
\end{displaymath}
The group of autotopisms of $L$ is denoted by $Aut(L, \cdot )$. If
$U\in S(L,\cdot )$ such that $(U,U,U)\in Aut(L,\cdot )$, then $U$ is
called an automorphism. Let $(L, \cdot )$ and $(G, \circ )$ be two
distinct loops. The triple
\begin{displaymath}
(U, V, W)~ :~ (L, \cdot )\to (G, \circ )~\textrm{such that}~U, V, W
: L\to G
\end{displaymath}
are bijections is called a loop isotopism if and only if
\begin{displaymath}
xU\circ yV=(x\cdot y)W~\forall ~x, y\in L.
\end{displaymath}
The representation sets of central loops are investigated and the
results obtained are used to construct a finite C-loop. We show that
for certain types of isotopisms, LC-identity and RC-identity are
isotopic invariant. For C-identity, this is true for commutative
loops in general. If the loops are non-commutative, then they must
be alternative central square for the C-identity to be an invariant
property. Under these types of isotopisms, alternative central
square loop isotopes of Moufang loops with 'lack of commutativity
property', some types of groups and RA-loops are found to be
C-loops.

\section{Preliminaries}

\begin{mydef}\label{definition:bijection}(\cite{phd3}, III.3.9~Definition, III.3.10~Definition, III.3.15~Definition)

Let $(L, \cdot )$ be a loop and $U, V, W\in S(L, \cdot )$.
\begin{enumerate}
\item If $(U, V, W)\in Aut(L, \cdot )$ for some $V, W$, then $U$ is called autotopic,
\begin{itemize}
\item the set of autotopic bijections in a loop $(L,\cdot )$ is represented by $\Sigma (L,\cdot )$.
\end{itemize}
\item If $(U, V, W)\in Aut(L, \cdot )$ such that $W=U, V=I$, then $U$ is called $\lambda $-regular,
\begin{itemize}
\item the set of all $\lambda $-regular bijections in a loop $(L,\cdot )$ is represented by $\Lambda (L,\cdot )$.
\end{itemize}
\item If $(U, V, W)\in Aut(L, \cdot )$ such that $U=I, W=V$, then $V$ is called $\rho $-regular,
\begin{itemize}
\item the set of all $\rho $-regular bijections in a loop $(L,\cdot )$ is represented by $P(L,\cdot )$.
\end{itemize}
\item If $\exists ~V\in S(L, \cdot )$ such that $xU\cdot y=x\cdot yV~\forall ~x, y\in L$,
then $U$ is called $\mu $-regular while $U'=V$ is called its adjoint.
\begin{itemize}
\item The set of all $\mu $-regular bijections in a loop $(L, \cdot )$ is denoted by $\Phi (L, \cdot )$,
while the collection of all adjoints in the loop $(L, \cdot )$ is denoted by $\Phi ^*(L, \cdot )$.
\end{itemize}
\end{enumerate}
\end{mydef}
The followings results will be judiciously used in this study.
\begin{myth}(\cite{phd3}, III.3.4~Theorem)

If two quasigroups are isotopic then their groups of autotopisms are isomorphic.
\end{myth}

\begin{myth}\label{lambdarhophi:subgroup}(\cite{phd3}, III.3.11~Theorem, III.3.16~Theorem)

The set $\Lambda (Q,\cdot )\Big(P(Q,\cdot )\Big)\Big[\Phi (Q,\cdot )\Big]$ of all $\lambda $-regular
($\rho $-regular)[$\mu $-regular] bijections of a quasigroup $(Q,\cdot )$ is a subgroup of the group
$\Sigma (Q,\cdot )$ of all autotopic bijections of $(Q,\cdot )$.
\end{myth}

\begin{mycor}\label{lambdarhophi:isomorphism}(\cite{phd3}, III.3.12~Corollary, III.3.16~Theorem)

If two quasigroups $Q$ and $Q'$ are isotopic, then the corresponding groups $\Lambda $
and $\Lambda '$($P$ and $P'$)[$\Phi $ and $\Phi '$]$\{\Phi ^*$ and $\Phi '^*\}$
are isomorphic.
\end{mycor}

\section{Main Results}
\subsection{The Representation sets of Central Loops}
\begin{myth}\label{lcrc:rep}
Let $\Pi_\lambda(\Pi_\rho)$ be the left(right) representation of a loop $L$.~$L$ is a LC(RC)-loop $\Leftrightarrow$
$\alpha,\beta\in\Pi_\lambda(\Pi_\rho)\Rightarrow \alpha\beta^2\in\Pi_\lambda(\Pi_\rho)$.
\end{myth}
{\bf Proof}\\ Let $L$ be an LC-loop, then $(x\cdot xy)z=x(x\cdot yz)$.
Using the definition of $L_x$, we have $L_{x\cdot xy}=L_yL_x^2$.
Replacing $x\cdot xy$ in $L$ and making $\alpha=L_y$ and $\beta=L_x$, $\alpha\beta^2\in\Pi_\lambda$.

Conversely, do the reverse of the above.
For $\Pi_\rho$ when $L$ is an RC-loop, $z(yx\cdot x)=(zy\cdot x)x$.
The proof goes in the same manner by using $R_x$.

\begin{myth}\label{c:rep}
Let $\Pi_\lambda(\Pi_\rho)$ be the left(right) representation of a loop $L$.~$L$ is a C-loop$\Leftrightarrow$
$\alpha,\beta\in\Pi_\lambda(\Pi_\rho)\Rightarrow \alpha\beta^2,\alpha^2\beta\in\Pi_\lambda(\Pi_\rho)$.
\end{myth}
{\bf Proof}\\ Let $L$ be a C-loop, then $(yx\cdot x)z=y(x\cdot
xz)\Rightarrow L_{yx\cdot x}=L_x^2L_y$. Let $\alpha=L_x ,
\beta=L_y$, then replacing $yx\cdot x$ in $L$,
$\alpha^2\beta\in\Pi_\lambda$. $L$ is a C-loop~$\Leftrightarrow L$
is an RC-loop and LC-loop by \cite{phd56}, hence by
Theorem~\ref{lcrc:rep}, $\alpha\beta^2,\alpha^2\beta\in\Pi_\lambda$.

Conversely let $\alpha,\beta\in\Pi_\lambda\Rightarrow \alpha\beta^2,\alpha^2\beta\in\Pi_\lambda$.
Take $\alpha=L_x,\beta=L_y$ then $(yx\cdot x)z=y(x\cdot xz)\Rightarrow~L$ is a C-loop.
For $\Pi_\rho$ the procedure is similar with $R_x$.

\begin{myth}\label{lcrc:power}
If $\Pi_\lambda(\Pi_\rho)$ is the left(right) representation of a LC(RC)-loop and
$\beta\in\Pi_\lambda(\Pi_\rho)$ then $\beta^n\in\Pi_\lambda(\Pi_\rho)~\forall~n\in \mathbb{Z}$.
\end{myth}
{\bf Proof}\\ By Theorem~\ref{lcrc:rep},
$\alpha,\beta\in\Pi_\lambda(\Pi_\rho)\Rightarrow
\alpha\beta^2\in\Pi_\lambda(\Pi_\rho)$. Using induction, when
$\alpha=I$, $\beta^2\in\Pi_\lambda(\Pi_\rho)$, when $\alpha=\beta$,
$\beta^3\in\Pi_\lambda(\Pi_\rho)$ and when $\alpha=\beta^k$,
$\beta^{k+2}\in\Pi_\lambda(\Pi_\rho)$ hence,
$\beta^n\in\Pi_\lambda(\Pi_\rho)~\forall~n\in \mathbb{Z}^+$. If $L$
is an LC(RC)-loop then it is a left(right) inverse property loop by
\cite{phd56}. Hence
$\beta=L_{x^{-1}}(R_{x^{-1}})\in\Pi_\lambda(\Pi_\rho)\Rightarrow\beta=L_x^{-1}(R_x^{-1})\in\Pi_\lambda(\Pi_\rho)~\forall~x\in
L$. By the earlier result,~$\beta^n\in\Pi_\lambda(\Pi_\rho)~\forall
~n\in \mathbb{Z}^-$. Whence
$\beta^n\in\Pi_\lambda(\Pi_\rho)~\forall~n\in \mathbb{Z}$.

\begin{mycor}\label{c:power}
If $\Pi_\lambda(\Pi_\rho)$ is the left(right) representation of a C-loop and
$\alpha\in\Pi_\lambda(\Pi_\rho)$ then $\alpha^n\in\Pi_\lambda(\Pi_\rho)~\forall~n\in \mathbb{Z}$.
\end{mycor}
{\bf Proof}\\ Using Theorem~\ref{c:rep}, $\alpha,\beta\in\Pi_\lambda(\Pi_\rho)\Rightarrow \alpha\beta^2,\alpha^2\beta\in\Pi_\lambda(\Pi_\rho)$.
The rest of the proof is similar to that in Theorem~\ref{lcrc:power}.

\subsection{Isotopes of LC, RC and C-loops}
\paragraph{}
Throughout this subsection, the following notations for translations
will be adopted; $L_x~:~y\mapsto xy$ and $R_x~:~y\mapsto yx$ for a
loop while $L_x'~:~y\mapsto xy$ and $R_x'~:~y\mapsto yx$ for its
loop isotope.

\begin{myth}\label{lc:auto}
A loop $L$ is an LC-loop $\Leftrightarrow (L_x^2,I,L_x^2)\in
Aut(L)~\forall~x\in L$.
\end{myth}
{\bf Proof}\\ Let $L$ be an LC-loop $\Leftrightarrow (x\cdot
xy)z=(xx)(yz)\Leftrightarrow(x\cdot xy)z=x(x\cdot yz)$ by \cite{den}
$\Leftrightarrow(L_x^2,I,L_x^2)\in Aut(L)~\forall~x\in L$.

\begin{myth}\label{rc:auto}
A loop $L$ is an RC-loop $\Leftrightarrow(I,R_x^2,R_x^2)\in
Aut(L)~\forall~x\in L$.
\end{myth}
{\bf Proof}\\ Let $L$ be an RC-loop, then $z(yx\cdot x)=zy\cdot
xx~\Leftrightarrow y(yx\cdot x)=(zy\cdot x)x$ by \cite{den}
$\Leftrightarrow(I,R_x^2,R_x^2)\in Aut(L)~\forall~x\in L$.

\begin{mylem}\label{lcrc:lp}
A loop is an LC(RC)-loop $\Leftrightarrow L_x^2(R_x^2)$ is $\lambda(\rho)$-regular
i.e $L_x^2(R_x^2)\in\Lambda (L)(P(L))$.
\end{mylem}
{\bf Proof}\\ Using Theorem~\ref{lc:auto}(Theorem~\ref{rc:auto}), the rest follows from the definition of $\lambda(\rho)$-regular bijection.

\begin{myth}\label{c:m}
A loop $L$ is a C-loop $\Leftrightarrow R_x^2$ is $\mu$-regular
and the adjoint of $R_x^2$, denoted by $(R_x^2)^*=L_x^2$ i.e $R_x^2\in\Phi (L)$ and $L_x^2\in\Phi ^*(L)$.
\end{myth}
{\bf Proof}\\ Let $L$ be a C-loop then $(yx\cdot x)z=y(x\cdot xz)\Rightarrow yR_x^2\cdot z=y\cdot zL_x^2\Rightarrow R_x^2\in \Phi (L)$ and $L_x^2\in \Phi ^*(L)$.
Conversely: do the reverse of the above.

\begin{myth}\label{iso:lcrc}
Let $(G, \cdot )$ and $(H, \circ )$ be any two distinct loops.
If the triple $\alpha =(A, B, B)~\Big(\alpha =(A, B, A)\Big)$ is an isotopism of $(G, \cdot )$ upon $(H, \circ )$,
then $(G, \cdot )$ is an LC(RC)-loop $\Leftrightarrow$ $(H, \circ )$ is a LC(RC)-loop.
\end{myth}
{\bf Proof}\\ By Lemma~\ref{lcrc:lp}, $G$ is an LC(RC)-loop $\Leftrightarrow L_x^2(R_x^2)\in \Lambda (G)(P(G))$.
Using the result in \cite{rit}, for each case ; $L_{xA}'=B^{-1}L_xB$ and $R_{xB}'=A^{-1}R_xA$~$\forall~x\in G$.
By Corollary~\ref{lambdarhophi:isomorphism}, there exists isomorphisms $\Lambda (G)\to \Lambda (H)$
and $P(G)\to P(H)$. Thus $L_y'^2(R_y'^2) \in \Lambda (H)(P(H))~\Leftrightarrow~H$ is an LC(RC)-loop.

\begin{myth}\label{iso:c}
Let $(G, \cdot )$ and $(H, \circ )$ be two distinct loops. If $G$ is a central square
C-loop, $H$ an alternative central square loop and the triple
$\alpha =(A, B, B)~\Big(\alpha =(A, B, A)\Big)$ is an isotopism of $G$ upon $H$, then
$H$ is a C-loop.
\end{myth}
{\bf Proof}\\ $G$ is a C-loop $\Leftrightarrow R_x^2\in\Phi (G)$ and $(R_x^2)^*=L_x^2\in\Phi ^*(G)$ for all $x\in G$
and using the result in \cite{rit} ;
$L_{xA}'=B^{-1}L_xB$ and $R_{xB}'=A^{-1}R_xA~\forall~x\in G$.
Using Corollary~\ref{lambdarhophi:isomorphism},
$R_y'^2 \in\Phi (H)$ and $(~R_y'^2~)^*=L_y'^2\in\Phi ^*(H) \Leftrightarrow H$ is a C-loop.

\begin{myth}\label{iso:cc}
Let $(G, \cdot )$ and $(H, \circ )$ be commutative loops.
If $\alpha =(A, B, B)$ or $\alpha =(A, B, A)$ is
an isotopism of $G$ upon $H$, then $G$ is a C-loop if and only if $H$ is a C-loop.
\end{myth}
{\bf Proof}\\ The proof is similar to that of Theorem~\ref{iso:c}

\begin{mycor}Under a triple of the form
$\alpha =(A,B,A)$ or $\alpha =(A,B,B)$, if a Dihedral group or
Quaternion group($Q_8$) or Cayley loop or $C_2\times C_2$ or
$C_2\times C_2\times C_2$ or RA-loops or $Q_8\times E\times A$ or
$M_{16}(Q_8)\times E\times A$ where $E$ is an elementary
abelian~2-group, $A$ is an abelian group(all of whose elements have
finite odd order) and $M_{16}(Q_8)$ is a Cayley loop, is isotopic to
an alternative central square loop $G$, then $G$ is a C-loop.
\end{mycor}
{\bf Proof}\\
From \cite{goo1} and \cite{goo2}, the Dihedral group $D_4$ of order
8, Quaternion group $Q_8$ of order 8, Cayley loop, $C_2\times C_2$,
$C_2\times C_2\times C_2$ and ring alternative loops(RA-loops) are
all central square. Hence, by Theorem~\ref{iso:c}, the claim that
$G$ is a C-loop follows. From \cite{nor} and the fact that $Q_8$ is
an Hamiltonian Moufang loop, $Q_8\times E\times A$ and
$M_{16}(Q_8)\times E\times A$ are both central square. Hence, by
Theorem~\ref{iso:c}, the claim that $G$ is a C-loop follows.

\paragraph{Construction}
Let $\Pi_\rho$ be the right representation of a
C-loop of order 12. If $\alpha, \beta, \gamma\in\Pi_\rho$ are given
by :
\begin{displaymath}
\alpha =(0~10~1~11~2~9)(3~7~4~8~5~6)=R_{10},
\end{displaymath}
\begin{displaymath}
\beta =(0~3)(1~4)(2~5)(6~10)(7~11)(8~9)=R_3,
\end{displaymath}
\begin{displaymath}
\gamma =(0~7~2~6~1~8)(3~10~5~9~4~11)=R_7.
\end{displaymath}
Then with Theorem~\ref{c:rep} and Corollary~\ref{c:power}, $\alpha,
\beta, \gamma\in\Pi_\rho$ generate other members of $\Pi_\rho$ by
considering the multiplications $\alpha\beta^2, \alpha^2\beta,
\gamma^n\in\Pi_\rho~\forall~n\in \mathbb{Z}$. These give us a finite
C-loop whose bordered multiplication table is shown by
Table~\ref{table:c} on Page~\pageref{table:c1}.
\begin{table}[!hbp]
\begin{center}
\begin{tabular}{|c||c|c|c|c|c|c|c|c|c|c|c|c|}
\hline
$\cdot $ & 0 & 1 & 2 & 3 & 4 & 5 & 6 & 7 & 8 & 9 & 10 & 11 \\
\hline  \hline
0 & 0 & 1 & 2 & 3 & 4 & 5 & 6 & 7 & 8 & 9 & 10 & 11 \\
\hline
1 & 1 & 2 & 0 & 4 & 5 & 3 & 7 & 8 & 6 & 10 & 11 & 9 \\
\hline
2 & 2 & 0 & 1 & 5 & 3 & 4 & 8 & 6 & 7 & 11 & 9 & 10 \\
\hline
3 & 3 & 4 & 5 & 0 & 1 & 2 & 9 & 10 & 11 & 6 & 7 & 8 \\
\hline
4 & 4 & 5 & 3 & 1 & 2 & 0 & 10 & 11 & 9 & 7 & 8 & 6 \\
\hline
5 & 5 & 3 & 4 & 2 & 0 & 1 & 11 & 9 & 10 & 8 & 6 & 7 \\
\hline
6 & 6 & 7 & 8 & 10 & 11 & 9 & 0 & 1 & 2 & 5 & 3 & 4 \\
\hline
7 & 7 & 8 & 6 & 11 & 9 & 10 & 1 & 2 & 0 & 3 & 4 & 5 \\
\hline
8 & 8 & 6 & 7 & 9 & 10 & 11 & 2 & 0 & 1 & 4 & 5 & 3 \\
\hline
9 & 9 & 10 & 11 & 8 & 6 & 7 & 3 & 4 & 5 & 2 & 0 & 1 \\
\hline
10 & 10 & 11 & 9 & 6 & 7 & 8 & 4 & 5 & 3 & 0 & 1 & 2 \\
\hline
11 & 11 & 9 & 10 & 7 & 8 & 6 & 5 & 3 & 4 & 1 & 2 & 0 \\
\hline
\end{tabular}
\end{center}\label{table:c1}
\caption{A non-associative C-loop of order 12}\label{table:c}
\end{table}

\end{document}